\documentclass[a4paper,12pt]{article}

\usepackage{amsfonts}
\usepackage{latexsym}
\usepackage{epsfig}
\usepackage{amssymb}
\usepackage{graphicx}
\usepackage{oldgerm}
\usepackage{theorem}

\def\RV{{\bf R}^{N}_{<}}

\def\a{{\bf a}}

\def\x{{\bf x}}
\def\y{{\bf y}}
\def\z{{\bf z}}

\def\cN{{\cal N}_N}
\def\H{{\cal H}(N)}
\def\S{{\cal S}(N)}
\def\U{{\bf U}(N)}
\def\O{{\bf O}(N)}
\def\dH{{\cal U}(dH)}
\def\dS{{\cal V}(dA)}

\def\R{{\bf R}}

\def\rR{{\rm R}}
\def\rI{{\rm I}}

\def\vlambda{\mib{\lambda}}

\theorembodyfont{\itshape}

\newtheorem{thm}{Theorem}[section]
\newtheorem{lem}[thm]{Lemma}
\newtheorem{cor}[thm]{Corollary}
\newtheorem{prop}[thm]{Proposition}

\newcommand{\mib}[1]{\mbox{\boldmath $#1$}}
\newcommand{\SSC}[1]{\section{#1}\setcounter{equation}{0}}
\newcommand{\qed}{\hbox{\rule[-2pt]{3pt}{6pt}}}

\begin{document}
\begin{center}
{\bf \Large{ Noncolliding Brownian motions and
\\Harish-Chandra formula}}
\vskip 3mm

\noindent
{\sc Makoto Katori}\\
Department of Physics,\\
Faculty of Science and Engineering,\\
Chuo University, \\
Kasuga, Bunkyo-ku, \\
Tokyo 112-8551, Japan \\
e-mail: katori@phys.chuo-u.ac.jp 

\vskip 3mm
\noindent
{\sc Hideki Tanemura}\\
Department of Mathematics and Informatics,\\
Faculty of Science,\\
Chiba University, \\
1-33 Yayoi-cho, Inage-ku,\\
Chiba 263-8522, Japan\\
e-mail: tanemura@math.s.chiba-u.ac.jp\\

\end{center}

\pagestyle{plain}
\vskip 0.3cm

\noindent  
{\bf Abstract.} 
We consider a system of noncolliding Brownian motions
introduced in our previous paper, 
in which the noncolliding condition is imposed 
in a finite time interval $(0,T]$.
This is a temporally inhomogeneous diffusion process
whose transition probability density depends on a value of $T$,
and in the limit $T \to \infty$ 
it converges to a temporally homogeneous diffusion process
called Dyson's model of Brownian motions.
It is known that the distribution of particle positions in Dyson's model 
coincides with that of eigenvalues of a Hermitian matrix-valued
process, whose entries are independent Brownian motions.
In the present paper we construct such a Hermitian matrix-valued process,
whose entries are sums of
Brownian motions and Brownian bridges given 
independently of each other, 
that its eigenvalues are identically distributed
with the particle positions of 
our temporally inhomogeneous system of noncolliding Brownian motions.
As a corollary of this identification we derive the
Harish-Chandra formula for an integral over the unitary group.

\vskip 0.3cm

\vspace{1cm}

\noindent
{\bf Abbreviated title.}
{\it  }  



\vspace{3mm}
\normalsize
%
\SSC{Introduction}

Dyson introduced a Hermitian matrix-valued process
whose $ij$-entry equals to
$B_{ij}(t)/\sqrt{2} + \sqrt{-1}\widehat{B}_{ij}(t)/\sqrt{2}$,
if $1\leq i <j\leq N$,
and equals to $B_{ii}(t)$, if $i=j$,
where $B_{ij}(t), \widehat{B}_{ij}(t)$,
$1\le i \le j \le N$, are independent Brownian motions \cite{Dys62}.
He found that its eigenvalues perform the
Brownian motions with the drift terms
acting as the repulsive two-body forces proportional to 
the inverse of distances between them,
which is now called Dyson's model of Brownian motions.
A number of processes of eigenvalues have been 
studied for random matrices by Bru \cite{Bru89, Bru91}, 
Grabiner \cite{Gra99}, 
K\"{o}nig and O'Connell \cite{KO01}, and others, but
all of them are temporally homogeneous diffusion processes.
In the present paper we introduce a Hermitian matrix-valued process,
whose eigenvalues give a temporally inhomogeneous diffusion process.

Let ${\bf Y}(t)=(Y_1(t),Y_2(t),\dots,Y_N(t))$ be
the system of $N$ independent Brownian motions in $\R$
conditioned never to collide to each other. 
It is constructed by the $h$-transform, in the sense of Doob \cite{Do84},
of the absorbing Brownian motion in a Weyl chamber 
of type $A_{N-1}$,
\begin{equation}
\RV=\{\x\in {\bf R}^N : x_1 < x_2 < \dots < x_N \}
\label{eqn:RV}
\end{equation}
with its harmonic function 
\begin{equation}
h_N(\x) = \prod_{1\le i <j \le N}( x_j - x_i),
\label{eqn:hN}
\end{equation}
$\x \in \RV$. We can prove that it is identically 
distributed with Dyson's model of Brownian motion.
In our previous papers \cite{KT02a, KT02b},
we introduce another system of noncolliding Brownian motions 
${\bf X}(t)=(X_1(t),X_2(t),\dots,X_N(t))$,
in which Brownian motions do not collide with each other
in a finite time interval $(0,T]$.
This is a temporally inhomogeneous diffusion process, 
whose transition probability density depends on the value of $T$.
It is easy to see that it converges to the process ${\bf Y}(t)$
in the limit $T\to\infty$. Moreover, it was shown that
$P({\bf X}(\cdot)\in dw)$ is absolutely continuous with respect to 
$P({\bf Y}(\cdot)\in dw)$ and that, 
in the case ${\bf X}(0)={\bf Y}(0)={\bf 0}$,
the Radon-Nikodym density is given by a constant 
multiple of $1/h_N(w(T))$.
Since this fact can be regard as an $N$-dimensional generalization of 
the relation proved by Imhof \cite{Imh84}
between a Brownian meander, which is temporally
inhomogeneous, and a three-dimensional Bessel process,
we called it {\it generalized Imhof's relation} \cite{KT02b}.

The problem we consider in the present paper is 
to determine a matrix-valued process 
that realizes ${\bf X}(t)$ as the process of its eigenvalues.
We found a hint in Yor \cite{Yor92} to solve this problem :
equivalence in distribution between the square of
the Brownian meander and the sum of the squares of a two-dimensional
Bessel process and of an independent Brownian bridge.
We prepare independent Brownian bridges 
$\beta_{ij}(t)$, $1\le i \le j \le N$ of duration $T$,
which are independent of the Brownian motions
$B_{ij}(t)$, $1\le i \le j \le N$, and set
a Hermitian matrix-valued process $\Xi^T(t)$, $t\in [0,T]$,
such that its $ij$-entry equals to
$B_{ij}(t)/\sqrt{2}+ \sqrt{-1}\beta_{ij}(t)/\sqrt{2}$,
if $1\leq i <j\leq N$, and it equals to $B_{ii}(t)$, if $i=j$.
Then we can prove that the eigenvalues of the matrix $\Xi^T(t)$ realize
${\bf X}(t), t \in [0,T]$ (Theorem 2.2).
This result demonstrates the fact that a temporally inhomogeneous diffusion
process ${\bf X}(t)$ in the $N$ dimensional space
can be represented as a projection of a combination of
$N(N+1)/2$ independent Brownian motions and 
$N(N-1)/2$ independent Brownian bridges.

It is known that Brownian motions $B_{ij}(t), 1\leq i,j \leq N$
are decomposed orthogonally into
the Brownian bridges $B_{ij}(t) - (t/T)B_{ij}(T)$
and the processes $(t/T)B_{ij}(T)$
(see, for instance, \cite{RY98, Yor92}).
Then the process $\Xi^T(t)$ can be decomposed into two independent
matrix-valued processes $\Theta^{(1)}(t)$ and $\Theta^{(2)}(t)$
such that, for each $t$,
the former realizes the distribution of 
Gaussian unitary ensemble (GUE) of complex Hermitian 
random matrices and the latter does of the 
Gaussian orthogonal ensemble (GOE) of real symmetric 
random matrices, respectively.
This implies that the process $\Xi^T(t)$ is identified with
a two-matrix model studied by Pandey and Mehta \cite{MP83, PM83},
which is a one-parameter interpolation of GUE and GOE,
if the parameter of the model is appropriately 
related with time $t$.
In \cite{KT02a} we showed this equivalence by using 
the Harish-Chandra formula for an integral
over the unitary group \cite{HC57}.
The proof of Theorem 2.2 makes effective use of 
our generalized version of Imhof's relation 
and this equivalence is established. 
The Harish-Chandra formula is then derived as a corollary 
of our theorem (Corollary 2.3).

As clarified by this paper, the Harish-Chandra integral formula implies 
the equivalence between temporally inhomogeneous systems of Brownian particles
and multi-matrix models. This equivalence is very useful
to calculate time-correlation functions 
of the particle systems. By using the method of orthogonal
polynomials developed in the random matrix theory \cite{Meh91},
determinantal expressions are derived for the correlations
and by studying their asymptotic behaviors, infinite
particle limits can be determined as reported in \cite{NKT03, KNT03}.

Extensions of the present results for the systematic study of relations
between noncolliding Brownian motions with geometrical restrictions
({\it e.g.} with an absorbing wall at the origin \cite{KT02b,KTNK03}) and
other random matrix ensembles than GUE and GOE
(see \cite{Meh91,Z96,AZ97}, for instance), 
will be reported elsewhere \cite{KT03b}.

\SSC{Preliminaries and Statement of Results}

\subsection{Noncolliding Brownian motions}

We consider the Weyl chamber of type $A_{N-1}$ as
(\ref{eqn:RV}) \cite{FH91, Gra99}.
By virtue of the Karlin-McGregor formula \cite{KM59_1,KM59_2}, 
the transition density function  $f_N (t, \y|\x)$ of
the absorbing Brownian motion in $\RV$
and the probability $\cN (t, \x)$ that the Brownian motion 
started at $\x\in\RV$ does not
hit the boundary of $\RV$ up to time $t>0$
are given by
\begin{equation}
\label{eqn:fnt}
f_{N}(t, \y|\x)= 
\det_{1 \leq i, j \leq N}
\Big[  G_t(x_j, y_i) \Big],
\: \x,\y \in \RV,
\end{equation}
and
$$
\label{eqn:cNt}
\cN (t, \x ) = 
\int_{\RV} d\y f_N (t, \y |\x ),
$$
respectively,
where $G_t(x,y)=(2 \pi t)^{-1/2} \ e^{-(y-x)^2/2t}$.
The function $h_N(\x )$ given by (\ref{eqn:hN})
is a strictly positive harmonic function for
absorbing Brownian motion in the Weyl chamber.
We denote by ${\bf Y}(t)=(Y_1(t),Y_2(t),\dots,Y_N(t)), t\in [0,\infty)$
the corresponding Doob $h$-transform \cite{Do84}, that is
the temporally homogeneous diffusion process
with transition probability density $p_N(s,\x,t,\y)$;
\begin{equation}
p_{N}(0, {\bf 0}, t, \y)
= \frac{t^{-N^2/2}}{ C_1(N)}
\exp\left\{ -\frac{|\y|^2}{2t} \right\} h_N(\y)^2,
\label{eqn:pn0}
\end{equation}
\begin{equation}
p_N(s,\x, t, \y )
= \frac{1}{h_N(\x)}f_{N}(t-s, \y|\x)h_N(\y),
\label{eqn:pnst}
\end{equation}
for $0 < s < t < \infty,\; \x, \y \in \RV,$
where $C_1(N) =(2\pi)^{N/2}\prod_{j=1}^{N} \Gamma(j)$.
The process ${\bf Y}(t)$ represents the system of $N$ Brownian motions
conditioned never to collide.
The diffusion process ${\bf Y}(t)$ 
solves the equation of Dyson's Brownian motion model \cite{Dys62} :
\begin{equation}
\label{eq:SDE2}
dY_i (t) = dB_i(t) 
+ \sum_{1\le j \le N,j\not=i} \frac{1}{Y_i(t)-Y_j(t)}dt,
\quad t\in [0,\infty), \quad i=1,2,\dots, N,
\end{equation}
where $B_i(t)$, $i=1,2,\dots,N$, are independent one dimensional 
Brownian motions.

For a given $T>0$, we define 
\begin{equation}
\label{eqn:gnst}
g_N^T(s,\x, t, \y )
= \frac{f_{N}(t-s, \y|\x)\cN (T-t,\y)}{\cN (T-s,\x)},
\end{equation}
for $ 0 < s < t \le T,\; \x, \y \in \RV$.
The function $g_{N}^{T}(s,\x,t,\y)$ can be
regarded as the transition probability density from the
state $\x \in \RV $ at time $s$ 
to the state $\y \in \RV$ at time $t$, 
and associated with the temporally inhomogeneous diffusion process, 
${\bf X}(t)=(X_1(t),X_2(t),\dots,X_N(t)), t\in [0,T]$,
which represents the system of $N$ Brownian motions conditioned not 
to collide with each other in a finite time interval $[0,T]$.
It was shown in \cite{KT02b} that
as $|\x|\to 0$, $g_{N}^{T}(0,\x,t,\y)$ converges to 
\begin{equation}
\label{eqn:gn0}
g_{N}^{T}(0, {\bf 0}, t, \y)
=\frac{T^{N(N-1)/4}t^{-N^2/2}}{C_2(N)}
\exp\left\{ -\frac{|\y|^2}{2t} \right\}
h_N(\y)\cN (T-t,\y),
\end{equation}
where 
$C_2(N) = 2^{N/2}\prod_{j=1}^N \Gamma(j/2)$.
Then the diffusion process ${\bf X}(t)$ solves the following equation:
\begin{equation}
\label{eqn:inhom1}
dX_i (t) = dB_i(t) + b_i^T(t, {\bf X}(t)) dt,
\quad t\in [0,T],\quad i=1,2,\dots, N,
\end{equation}
where 
$$
b_i^T (t,\x) = \frac{\partial}{\partial x_i}\ln \cN (T-t, \x),
\qquad i=1,2,\dots,N.
$$

From the transition probability densities 
(\ref{eqn:pn0}), (\ref{eqn:pnst}) 
and (\ref{eqn:gn0}), (\ref{eqn:gnst})
of the processes,
we have the following relation between 
the processes ${\bf X}(t)$ and ${\bf Y}(t)$ 
in the case ${\bf X}(0)={\bf Y}(0)={\bf 0}$ \cite{KT02a, KT02b}:

\begin{equation}
\label{eqn:Imhof}
P( {\bf X}(\cdot) \in dw)
= \frac{ C_1(N) T^{N(N-1)/4}}{C_2(N) h_N (w(T))}
P( {\bf Y}(\cdot) \in dw).
\end{equation}
This is the generalized form of 
the relation obtained by Imhof \cite{Imh84} for
the Brownian meander and the three-dimensional Bessel process.
Then, we call it {\it generalized Imhof's relation}.
\vskip 3mm

\subsection{Hermitian matrix-valued processes}

We denote by $\H$ the set of $N\times N$ complex Hermitian matrices 
and by $\S$ the set of $N\times N$ real symmetric matrices.
We consider complex-valued processes $x_{ij}(t)$, $1\le i,j \le M$
with $x_{ij}(t)= x_{ji}(t)^{\dagger}$,
and Hermitian matrix-valued processes
$\Xi (t)= (x_{ij}(t))_{1\le i,j\le N}$.

Here we give two examples of Hermitian matrix-valued process.
Let $B^{\rR}_{ij}(t)$, $B^{\rI}_{ij}(t)$, $1\le i \leq j \leq N$
be independent one dimensional Brownian motions. 
Put
$$
x_{ij}^{\rR}(t)
=
\left\{
   \begin{array}{ll}
      \displaystyle{
      \frac{1}{\sqrt{2}} B_{ij}^{\rR}(t),
      } & 
   \mbox{if} \ i < j, \\
        & \nonumber\\
        B_{ii}^{\rR}(t), & \mbox{if} \ i=j,
     \nonumber\\
   \end{array}\right. 
   \quad {\rm and } \quad
x_{ij}^{\rI}(t)
=
\left\{
   \begin{array}{ll}
      \displaystyle{
      \frac{1}{\sqrt{2}} B_{ij}^{\rI}(t),
      } & 
   \mbox{if} \ i < j, \\
        & \nonumber\\
    0, & \mbox{if} \ i=j,
     \nonumber\\
   \end{array}\right. 
$$
with $x_{ij}^{\rR}(t)=x_{ji}^{\rR}(t)$ 
and $x_{ij}^{\rI}(t)=- x_{ji}^{\rI}(t)$ for $i > j$.


\vskip 3mm

\noindent (i) {\bf GUE type matrix-valued process.}
\ Let 
$\Xi^{\rm GUE} (t) 
=(x_{ij}^{\rR}(t)+\sqrt{-1}x_{ij}^{\rI}(t))_{1 \leq i, j \leq N}$.
For fixed $t\in [0,\infty)$, $\Xi^{\rm GUE} (t)$
is in the Gaussian unitary ensemble (GUE), that is,
its probability density function with respect to 
the volume element $\dH$ of $\H$ is given by
$$
\mu^{\rm GUE}(H , t)
= 
\frac{t^{-N^2/2}}{C_3(N)}
\exp \left\{ - \frac{1}{2t} {\rm Tr}H^2 \right\},
\quad H\in \H,
$$
where $C_3(N)= 2^{N/2} \pi^{N^2/2}$.
Let $\U$ be the space of all $N\times N$ unitary matrices.
For any $U\in\U$,
the probability $\mu^{\rm GUE}(H, t)\dH $ is invariant under 
the automorphism $H\to U^{\dagger}HU$.
It is known that the distribution of eigenvalues
$\x \in \RV$ of the matrix ensembles is given as
$$
g^{\rm GUE}(\x , t)
= \frac{t^{-N^2/2}}{C_1(N)} 
\exp \left\{ - \frac{|\x|^2}{2t}\right\}
h_N(\x)^2,
$$
\cite{Meh91}, and so 
$p_N(0,{\bf 0},t,\x)=g^{\rm GUE}(\x , t)$ from (\ref{eqn:pn0}).

\vskip 3mm

\noindent (ii) {\bf GOE type matrix-valued process.}
\ Let $\Xi^{\rm GOE} (t) =(x_{ij}^{\rR}(t))_{1 \leq i, j \leq N}$.
For fixed $t\in [0,\infty)$, $\Xi^{\rm GOE} (t)$
is in the Gaussian orthogonal ensemble (GOE),
that is, 
its probability density function with respect to 
the volume element $\dS$ of $\S$ is given by
$$
\mu^{\rm GOE}(A , t)
= 
\frac{t^{-N(N+1)/4}}{C_4(N)}
\exp \left\{ - \frac{1}{2t} {\rm Tr}A^2 \right\},
\quad A\in \S,
$$
where $C_4(N)= 2^{N/2} \pi^{N(N+1)/4}$.
Let $\O$ be the space of all $N\times N$ real orthogonal matrices.
For any $V\in\O$,
the probability $\mu^{\rm GOE}(H, t)\dS$ is invariant under 
the automorphism $A \to V^{T}AV$.
It is known that the probability density of eigenvalues
$\x \in \RV$ of the matrix ensemble is given as
$$
g^{\rm GOE}(\x , t)
= \frac{t^{-N(N+1)/4}}{C_2(N)} 
\exp \left\{ - \frac{|\x|^2}{2t}\right\}
h_N(\x),
$$
\cite{Meh91}, and so $g_N^t(0,{\bf 0},t,\x)=g^{\rm GOE}(\x , t)$
from (\ref{eqn:gn0}).

\vskip 3mm

We denote by $U(t)=(u_{ij}(t))_{1\le i,j\le N}$ 
the family of unitary matrices which diagonalize $\Xi (t)$:
$$
U(t)^{\dagger} \Xi (t) U(t) = \Lambda (t) = {\rm diag}\{ \lambda_i(t) \},
$$
where $\{\lambda_i(t)\}$ are eigenvalues of $\Xi (t)$
such that $\lambda_1(t) \leq \lambda_2(t) \leq \cdots \leq \lambda_N(t)$.
By a slight modification of Theorem 1 in Bru \cite{Bru89}
we have the following.

\vskip 3mm


\begin{prop}
\label{prop:Bru}
Let $x_{ij}(t)$, $1\leq i,j \leq N$ be continuous semimartingales.
The process $\vlambda (t) = 
(\lambda_1(t),\lambda_2(t),\dots, \lambda_N(t))$
satisfies
\begin{equation}
\label{eqn:l=M+J}
d\lambda_i(t) = dM_i(t) + dJ_i(t), 
\quad i=1,2,\dots, N,
\end{equation}
where $M_i(t)$ is the martingale with quadratic variation 
$\langle M_i \rangle_t = \int_0^t \Gamma_{ii}(s)ds$,
and $J_i(t)$ is the process with finite variation given by
\begin{eqnarray}
d J_i(t)&=& 
\sum_{j=1}^{N} \frac{1}{\lambda_i(t) -\lambda_j(t)}
1(\lambda_i \not= \lambda_j) \Gamma_{ij}(t)dt
\nonumber
\\
&+& {\it the \ finite \ variation \ part \ of }\ 
(U(t)^{\dagger} d \Xi (t) U(t))_{ii}
\nonumber
\end{eqnarray}
with
\begin{equation}
\Gamma_{ij}(t)dt = 
(U^{\dagger}(t) d \Xi (t) U(t))_{ij}
(U^{\dagger}(t) d \Xi (t) U(t))_{ji}.
\end{equation}
\end{prop}

\vskip 3mm

For the process $\Xi^{\rm GUE}(t)$, 
$d\Xi^{\rm GUE}_{ij}(t)d\Xi^{\rm GUE}_{k \ell}(t)
= \delta_{i \ell} \delta_{jk} dt$
and $\Gamma_{ij}(t) =1$.
The equation (\ref{eqn:l=M+J}) is given as
$$
d \lambda_{i}(t)=d B_{i}(t)+
\sum_{j : j \not=i} \frac{1}{\lambda_{i}(t)-\lambda_{j}(t)}dt,
\quad 1\leq i \leq N.
$$
Hence, the process $\lambda (t)$ is the homogeneous diffusion
that coincides with the system of noncolliding Brownian motions 
${\bf Y}(t)$ with ${\bf Y}(0)={\bf 0}$.


For the process $\Xi^{\rm GOE}(t)$, 
$d\Xi^{\rm GOE}_{ij}(t)d\Xi^{\rm GOE}_{k \ell}(t)
= \frac{1}{2}\Big( \delta_{i \ell} \delta_{jk} 
+\delta_{ik} \delta_{j \ell} \Big) dt$
and 
$\Gamma_{ij}(t)= \frac{1}{2}(1+\delta_{ij})$.
The equation (\ref{eqn:l=M+J}) is given as
$$
d \lambda_{i}(t)=d B_{i}(t)+ \frac{1}{2}
\sum_{j : j \not=i} \frac{1}{\lambda_{i}(t)-\lambda_{j}(t)}dt,
\quad 1\leq i \leq N.
$$

\subsection{Results}

Let $\beta_{ij}(t)$, $1\leq i < j \leq N$ 
be independent one dimensional Brownian bridges of duration $T$,
which are the solutions of the following equation:
$$
\beta_{ij}(t) = B^{\rI}_{ij}(t) - \int_0^t \frac{\beta_{ij}(s)}{T-s}ds,
\quad 0\leq t \leq T. 
$$
For $t\in [0,T]$, we put 
$$
\xi_{ij}(t)
=
\left\{
   \begin{array}{ll}
      \displaystyle{
      \frac{1}{\sqrt{2}}\beta_{ij}(t),
      } & 
   \mbox{if} \ i < j, \\
        & \nonumber\\
      0, & \mbox{if} \ i=j,
     \nonumber\\
   \end{array}\right. 
$$
with $\xi_{ij}(t) = - \xi_{ji}(t)$ for $i>j$.
We introduce the $\H$-valued process
$\Xi^{T}(t)= (x_{ij}^{\rR}(t)+\sqrt{-1}\xi_{ij}(t))_{1\le i,j\le N}$.
Then, the main result of this paper is the following theorem.


\begin{thm}
\label{thm:th_inhom}
Let $\lambda_i(t)$, $i=1,2,\dots,N$ be the eigenvalues of $\Xi^{T}(t)$
with $\lambda_1(t) \leq \lambda_2(t) \leq \cdots \leq \lambda_N(t)$.
The process $\vlambda(t) = 
(\lambda_1(t),\lambda_2(t), \dots,\lambda_N(t))$ 
is the temporally inhomogeneous diffusion that coincides with
the noncolliding Brownian motions ${\bf X}(t)$ with ${\bf X}(0)={\bf 0}$.
\end{thm}

As a corollary of the above result,
we have the following formula,
which is called the Harish-Chandra
integral formula \cite{HC57}
(see also \cite{IZ80,Meh91}).
Let $dU$ be the Haar measure of the space ${\bf U}(N)$
normalized as $\int_{\U}dU =1$.

\vskip 0.5cm

\begin{cor}
\label{thm:HCIZ}
Let $\x =(x_1,x_2,\dots,x_N), \y=(y_1,y_2,\dots,y_N) \in \RV$.
Then
$$
 \int_{\U} dU \, \exp  \left\{ -\frac{1}{2\sigma^2} {\rm Tr}
(\Lambda_{\x}-U^{\dagger} \Lambda_{\y} U)^2 \right\}
\nonumber\\
= \frac{C_1(N)\sigma^{N^2}}{h_{N}(\x) h_{N}(\y)}
\det_{1 \leq i, j \leq N} 
\Big[ G_{\sigma^2} (x_{i}, y_{j}) \Big],
$$
where 
$\Lambda_{\x}={\rm diag}\{ x_{1}, \dots, x_{N}\}$ and
$\Lambda_{\y}={\rm diag}\{ y_{1}, \dots, y_{N}\}$.
\end{cor}

\noindent
{\bf Remark} \
Applying Proposition \ref{prop:Bru}
we derive the following equation:
\begin{equation}
\label{eqn:inhom2}
d\lambda_{i}(t) = dB_i(t) 
+ \sum_{j: j\not= i}\frac{1}{\lambda_i(t)-\lambda_j(t)}dt
- \frac{\lambda_i(t)-
\int_{\S} \mu^{\rm GOE}(dA)(U(t)^{\dagger} A U(t))_{ii}}{T-t} dt,
\end{equation}
$i=1,2,\dots,N$, where $U(t)$ is one of the families of 
unitary matrices which diagonalize $\Xi^T (t)$.
From the equations (\ref{eqn:inhom1}) and (\ref{eqn:inhom2})
we have
\begin{eqnarray}
&&\int_{\S} \mu^{\rm GOE}(dA)(U(t)^{\dagger} A U(t))_{ii}
\nonumber\\
&&= \lambda_i(t) +(T-t) \left\{
\frac{\frac{\partial}{\partial \lambda_i}
\cN (T-t, \vlambda (t))}{\cN (T-t, \vlambda (t))}
-\sum_{j: j\not = i}\frac{1}{\lambda_i(t) -\lambda_j(t)}
\right\}.
\label{eqn:drift}
\end{eqnarray}
The function $\cN (t, \x)$ is expressed by a Pfaffian of the matrix
whose $ij$-entry is $\Psi ((x_j-x_i)/2\sqrt{t})$
with $\Psi (u) = \int_0^u e^{-v^2} dv$.
(See Lemma 2.1 in \cite{KT02b}.)
Then the right hand side of (\ref{eqn:drift}) can be written explicitly.

\SSC{Proofs}

\subsection{Proof of Theorem \ref{thm:th_inhom}}

For $y\in\R$ and $1\le i,j \le N$,
let $\beta_{ij}^{\sharp}(t)=\beta_{ij}^{\sharp}(t:y)$, $t\in [0,T]$,
$\sharp =\rR,\rI$, be diffusion processes 
which satisfy the following stochastic differential equations:
\begin{equation}
\beta_{ij}^{\sharp}(t : y) 
= B_{ij}^{\sharp}(t) 
- \int_0^t \frac{\beta^{\sharp}_{ij}(s : y) - y}{T-s}ds,
\quad t\in [0,T].
\label{eqn:bbridge_1}
\end{equation}
These processes are Brownian bridges of duration $T$ 
starting form $0$ and ending at $y$.
For $H = (y_{ij}^{\rR} + \sqrt{-1}y_{ij}^{\rI})_{1\leq i,j \leq N}\in \H$ 
we put 
\begin{eqnarray}
&&\xi_{ij}^{\rR}(t:y_{ij}^{\rR}) 
= 
\left\{
   \begin{array}{ll}
      \displaystyle{
      \frac{1}{\sqrt{2}}\beta^{\rR}_{ij}(t:\sqrt{2}y_{ij}^{\rR}),
      } & 
   \mbox{if} \ i < j, \\
        & \nonumber\\
\beta^{R}_{ii}(t:y_{ii}^{\rR}),& \mbox{if} \ i=j,
     \nonumber\\
   \end{array}\right.
\\
&&\xi_{ij}^{\rI}(t:y_{ij}^{\rI}) 
= 
\left\{
   \begin{array}{ll}
      \displaystyle{
      \frac{1}{\sqrt{2}}\beta_{ij}^{\rI}(t:\sqrt{2}y_{ij}^{\rI}),
      } & 
   \mbox{if} \ i < j, \\
        & \nonumber\\
0,& \mbox{if} \ i=j,
     \nonumber\\
\end{array}\right.
\end{eqnarray}
with $\xi_{ij}^{\rR}(t:y_{ij}^{\rR})= \xi_{ji}^{\rR}(t:y_{ji}^{\rR})$
and 
$\xi_{ij}^{\rI}(t:y_{ij}^{\rI})=- \xi_{ji}^{\rI}(t:y_{ji}^{\rI})$ for $i>j$.
We introduce the  $\H$-valued process
$\Xi^T(t:H) = ( \xi_{ij}^{\rR}(t:y_{ij}^{\rR}) 
+ \sqrt{-1} \xi_{ij}^{\rI}(t:y_{ij}^{\rI}) )_{1\leq i,j \leq N}$, 
$t\in [0,T]$.
From the equation (\ref{eqn:bbridge_1}) 
we have the equality
\begin{equation}
\label{eqn:bbridge_2}
\Xi^T (t:H) = \Xi^{\rm GUE}(t) - \int_0^t \frac{\Xi^T (s :H)-H}{T-s}ds,
\; t\in [0,T].
\end{equation}
Let $H_U$ be a random matrix with distribution $\mu^{\rm GUE}(\cdot,T)$,
and $A_O$ be a random matrix with distribution $\mu^{\rm GOE}(\cdot,T)$.
Note that $\beta^{\sharp}_{ij}(t:Y)$, $t\in [0,T]$ is a Brownian motion
when $Y$ is a Gaussian random variable with variance $T$,
which is independent of $B^{\sharp}_{ij}(t), t\in [0,T]$.
Then when $H_U$ and $A_O$ are independent of 
$\Xi^{\rm GUE} (t), t\in [0,T]$, 
\begin{eqnarray}
\label{eqn:H_U}
&&\Xi^T (t:H_U) = \Xi^{\rm GUE}(t),
\quad t \in [0,T],
\\
\label{eqn:A_O}
&&\Xi^T (t:A_O) = \Xi^T(t),
\quad t \in [0,T],
\end{eqnarray}
in the sense of distribution.
Since the distribution of the process $\Xi^{\rm GUE}(t)$ 
is invariant under any unitary transformation,
we obtain the following lemma from (\ref{eqn:bbridge_2}).

\begin{lem}
For any $U\in \U$ we have
$$
U^{\dagger}\Xi^T(t:H)U = \Xi^T(t: U^{\dagger}H U),
\quad t\in [0,T],
$$
in distribution.
\end{lem}
\vskip 1mm

From the above lemma it is obvious that
if $H^{(1)}$ and $H^{(2)}$ are $N\times N$ Hermitian matrices
having the same eigenvalues,
the processes of eigenvalues of $\Xi^T(t:H^{(1)}), t\in [0,T]$ 
and $\Xi(t:H^{(2)}), t\in [0,T]$ are identical in distribution.
For an $N\times N$ Hermitian matrix $H$ 
with eigenvalues $\{ a_i \}_{1\leq i \leq N}$,
we denote the probability distribution of the process of the eigenvalues of 
$\Xi^T(t:H)$ by $Q_{0,\a}^T(\cdot), t\in [0,T]$.
We also denote by $Q^{\rm GUE}(\cdot)$
the distribution of the process of eigenvalues of 
$\Xi^{\rm GUE}(t), t\in [0,T]$, and by $Q^T(\cdot)$ that of 
$\Xi^T(t), t\in [0,T]$. From the equalities (\ref{eqn:H_U}) 
and (\ref{eqn:A_O}) we have
\begin{eqnarray}
&&Q^{\rm GUE}(\cdot)
= \int_{\bf R_{<}^N} Q_{0,\a}^T(\cdot)g^{\rm GUE}(\a,T)d\a,
\nonumber
\\
&&Q^{T}(\cdot)= \int_{\bf R_{<}^N} Q_{0,\a}^T(\cdot)g^{\rm GOE}(\a,T)d\a.
\nonumber
\end{eqnarray}
Since $Q^{\rm GUE}(\cdot)$ is the distribution of 
the temporally homogeneous diffusion process ${\bf Y}(t)$
which describes noncolliding Brownian motions,
by our generalized Imhof's relation (\ref{eqn:Imhof}) we 
can conclude that $Q^{T}(\cdot)$ is the distribution of the 
temporally inhomogeneous diffusion process
${\bf X}(t)$ which describes our noncolliding Brownian motions.
\qed
\vskip 3mm

\subsection{Proof of Corollary \ref{thm:HCIZ}}
By (\ref{eqn:gn0}) we have
\begin{eqnarray}
&&g_{N}^{T}(0, {\bf 0}, t, \y)
 = \frac{1}{C_2(N)}T^{N(N-1)/4}t^{-N^2/2}  
\exp\Big\{ -\frac{|\y|^2}{2t} \Big\}
h_N(\y)
\nonumber\\
&&\qquad \times 
\int_{\RV} d\z \det_{1 \leq i, j \leq N}
\left[ \frac{1}{\sqrt{2 \pi (T-t)}} \ 
\exp \left\{ - \frac{(y_{j}-z_{i})^2}{2(T-t)} \right\} \right]
\nonumber
\\
&&\quad = \frac{1}{C_2(N)}T^{N(N-1)/4}t^{-N^2/2}  (2 \pi (T-t))^{-N/2}h_N(\y)
\nonumber
\\
&&\qquad \times 
\int_{\RV} d\z \det_{1 \leq i, j \leq N}
\left[
\exp \left\{ -\frac{y_j^2}{2t}
- \frac{(y_{j}-z_{i})^2}{2(T-t)} \right\}
\right]
\nonumber
\\
&&\quad= \frac{1}{C_2(N)}T^{N(N-1)/4}t^{-N^2/2}  (2 \pi (T-t))^{-N/2}h_N(\y)
\nonumber
\\
&&\qquad\times 
\int_{\RV} d\z 
\exp \left\{ -\frac{|\z|^2}{2T}\right\}
\det_{1 \leq i, j \leq N}
\left[
\exp \left\{ 
- \frac{T}{2t(T-t)}\left( y_{j}-\frac{t}{T}z_{i}\right)^2 \right\}
\right].
\nonumber
\end{eqnarray}
Setting $(t/T)z_i =a_i$, $i=1,2,\dots,N$, 
$t(T-t)/T = \sigma^2$ and $T/t^2 = \alpha$,
we have
\begin{eqnarray}
&&g_{N}^{T}(0, {\bf 0}, t, \y) = 
\frac{(2\pi)^{-N/2}}{C_2(N)}\sigma^{-N} \alpha^{N(N+1)/4}h_N(\y)
\nonumber
\\
&&\qquad\times 
\int_{\RV} d\a
\exp \left\{ -\frac{\alpha}{2}|\a|^2 \right\}
\det_{1 \leq i, j \leq N}
\left[
\exp \left\{ 
- \frac{1}{2\sigma^2}\left( y_{j}-a_{i}\right)^2 \right\}
\right].
\label{eqn:HC_1}
\end{eqnarray}

We write the transition probability density of the process $\Xi^T(t)$
by $q^T_N(s, H_1, t, H_2)$, $0\le s < t \le T$,
for $H_1, H_2 \in \H$.
Then by Theorem \ref{thm:th_inhom} and the fact that
$\dH = C_U(N) h_N(\y)^2 dU d\y$, with $C_U(N) = C_3(N) / C_1(N)$,
we have
\begin{equation}
g^T_N(0,{\bf 0},t, \y)
= C_U(N) h_N(\y)^2 
\int_{\U} dU \ q^T_N (0, O, t, U^{\dagger} \Lambda_{\y} U),
\label{eqn:HC_2}
\end{equation}
where $O$ is the zero matrix.
We introduce the $\H$-valued process
$\Theta^{(1)}(t)=(\theta^{(1)}_{ij}(t))_{1\leq i,j \leq N}$ and 
the $\S$-valued process
$\Theta^{(2)}(t)=(\theta^{(2)}_{ij}(t))_{1\leq i,j \leq N}$
which are defined by
\begin{eqnarray}
&&\theta^{(1)}_{ij}(t) = 
\left\{
   \begin{array}{ll}
      \displaystyle{
\frac{1}{\sqrt{2}}\left\{ B^{\rR}_{ij}(t) - \frac{t}{T}B^{\rR}_{ij}(T)\right\}
+ \frac{\sqrt{-1}}{\sqrt{2}}\beta_{ij}(t),
} & 
   \mbox{if} \ i < j, \\
        & \nonumber\\
        \displaystyle{
B^{\rR}_{ii}(t)- \frac{t}{T}B^{\rR}_{ii}(T),}
& \mbox{if} \ i=j,
     \nonumber\\
   \end{array}\right. \nonumber
\end{eqnarray}
and
\begin{eqnarray}
&&\theta^{(2)}_{ij}(t) = 
\left\{
   \begin{array}{ll}
      \displaystyle{
\frac{t}{\sqrt{2}T}B^{\rR}_{ij}(T),
} & 
   \mbox{if} \ i < j, \\
& \nonumber\\
\displaystyle{
\frac{t}{T}B^{\rR}_{ii}(T)},
& \mbox{if} \ i=j,
     \nonumber\\
   \end{array}\right. \nonumber
\end{eqnarray}
respectively.
Then $\Xi^T(t) = \Theta^{(1)}(t) + \Theta^{(2)}(t)$.
Note that $B^{\rR}_{ij}(t) - (t/T)B^{\rR}_{ij}(T)$
are Brownian bridges of duration $T$
which are independent of $(t/T)B^{\rR}_{ij}(T)$.
Hence $\Theta^{(1)}(t)$ is in the GUE
and $\Theta^{(2)}(t)$ is in the GOE
independent of $\Theta^{(1)}(t)$.
Since $E[ \theta_{ii}^{(1)}(t)^2]=\sigma^2$
and $E[ \theta_{ii}^{(2)}(t)^2]= 1/\alpha$,
the transition probability density $q_N^T (0,O,t,H)$ can be written by
\begin{eqnarray}
&&q_N^T (0,O,t,H) 
= \int_{\S} \dS \ 
\mu^{\rm GOE}\left( A, \frac{1}{\alpha}\right)\mu^{\rm GUE}(H-A,\sigma^2)
\nonumber\\
&&\quad = \frac{C_O(N) \sigma^{-N^2}\alpha^{N(N+1)/4}}{C_3(N) C_4(N)}
\int_{\RV} d\a \ h_N(\a)
\exp \left\{ -\frac{\alpha}{2} |\a|^2 
-\frac{1}{2\sigma^2} {\rm Tr} (H-\Lambda_{\a})^2 \right\},
\label{eqn:HC_3}
\end{eqnarray}
where we used the fact $\dS = C_O(N) h_N(\a) dV d\a$ with
the Haar measure $dV$ of the space $\O$ 
normalized as $\int_{\O}dV =1$, and
$C_O (N) = C_4(N) / C_2 (N)$.
Combining (\ref{eqn:HC_1}), (\ref{eqn:HC_2}) and (\ref{eqn:HC_3}) 
we have
\begin{eqnarray}
&&\frac{C_1(N) \sigma^{N^2-N} }{(2\pi)^{N/2} h_N(\y)}
\int_{\RV} d\a
\exp \left\{ -\frac{\alpha}{2}|\a|^2 \right\}
\det_{1 \leq i, j \leq N}
\left[
\exp \left\{ 
- \frac{1}{2\sigma^2}\left( y_{j}-a_{i}\right)^2 \right\}
\right]
\nonumber\\
&&\quad = \int_{\RV} d\a \,
h_N(\a) \exp \left\{ -\frac{\alpha}{2} |\a|^2 \right\}
\int_{\U} dU
\exp \left\{ -\frac{1}{2\sigma^2} {\rm Tr} 
(U^{\dagger}\Lambda_{\y} U -\Lambda_{\a})^2 \right\}.
\label{eqn:HC_4}
\end{eqnarray}
For each $\sigma >0$, (\ref{eqn:HC_4}) holds 
for any $\alpha >0$ and we have
\begin{eqnarray}
&& 
\frac{C_1(N) \sigma^{N^2} }{h_N(\y)h_N(\a)}
\det_{1 \leq i, j \leq N}
\left[ \frac{1}{\sqrt{2\pi \sigma^2}}
\exp \left\{ 
- \frac{1}{2\sigma^2}\left( y_{j}-a_{i}\right)^2 \right\}
\right]
\nonumber\\
&&\quad = \int_{\U} dU
\exp \left\{ -\frac{1}{2\sigma^2} {\rm Tr} 
(U^{\dagger}\Lambda_{\y} U -\Lambda_{\a})^2 \right\}.
\nonumber
\end{eqnarray}
This completes the proof.
\qed
\vskip 3mm

\footnotesize 


\clearpage

\end{document}